\documentclass[11pt]{amsart}
\usepackage[margin=1.15in]{geometry}
\usepackage{amsmath,amssymb,amsthm,mathtools}
\usepackage{booktabs}
\usepackage{microtype}
\usepackage[dvipsnames]{xcolor}
\usepackage[T1]{fontenc}
\usepackage{lmodern}
\usepackage[colorlinks=true,linkcolor=blue!50!black,citecolor=blue!50!black,urlcolor=blue!50!black]{hyperref}

\newtheorem{theorem}{Theorem}[section]
\newtheorem{proposition}[theorem]{Proposition}
\newtheorem{lemma}[theorem]{Lemma}
\newtheorem{corollary}[theorem]{Corollary}
\newtheorem{conjecture}[theorem]{Conjecture}
\newtheorem*{theoremA}{Theorem A}
\newtheorem*{theoremB}{Theorem B}
\theoremstyle{definition}
\newtheorem{definition}[theorem]{Definition}
\newtheorem{remark}[theorem]{Remark}
\newtheorem{example}[theorem]{Example}
\newtheorem{problem}[theorem]{Problem}

\newcommand{\A}{\mathcal{A}}
\newcommand{\C}{\mathbb{C}}
\newcommand{\PP}{\mathbb{P}}
\newcommand{\Der}{\operatorname{Der}}
\newcommand{\exps}{\operatorname{exp}}

\title[On Terao's conjecture for at most nine planes]{On Terao's conjecture for at most nine planes}

\author{Bartosz Jaros\l awski}
\address{Bartosz Jaros\l awski: Department of Mathematics, University of the National Education Commission Krakow, Podchor\k{a}\.zych 2, PL-30-084 Krak\'ow, Poland.}
\email{}


\begin{document}

\begin{abstract}
In this note we show that Terao's freeness conjecture holds for arrangements of at most nine planes in the complex projective $3$-space. Moreover, we show that every essential and irreducible free arrangement of at most nine planes in $\PP^{3}_{\C}$ is divisionally free, and we list all intersection lattices of such arrangements. Our approach is based on Ziegler's restriction theorem and a simple count involving the numbers of multiple lines. As a by-product, we show that Cynk--Szemberg octic arrangements are never free. Finally, motivated by a recent work of Dimca and Pokora, we present examples of weak Ziegler pairs of plane arrangements, including pairs of octic arrangements leading to rigid double octic Calabi--Yau threefolds.
\end{abstract}

\maketitle

\section{Introduction}

Let $\A\subset\PP^{n}_{\C}$ be an arrangement of $d$ hyperplanes defined by $Q(\A)=\ell_1\cdots\ell_d$, where the $\ell_i$ are linear forms in $S=\C[x_0,\ldots,x_n]$. Recall that $\A$ is \emph{free} if the module of logarithmic derivations
\[
D(\A)=\{\theta\in\Der(S)\,:\,\theta(Q(\A))\in Q(\A)\cdot S\}
\]
is a free $S$-module. Freeness is an algebraic property of $\A$, while the intersection lattice $L(\A)$ is purely combinatorial. The following conjecture, due to Terao, is one of the central open problems in the theory of hyperplane arrangements.

\begin{conjecture}[Terao]
Let $\A$ and $\mathcal{B}$ be two arrangements in $\PP^n_{\C}$ with $L(\A)\cong L(\mathcal{B})$. If $\A$ is free, then $\mathcal{B}$ is free.
\end{conjecture}

The conjecture is known to hold for line arrangements with a small number of lines. Faenzi and Vall\`es proved it for arrangements of up to $12$ lines \cite{FV}; Dimca, Ibadula and M\u{a}cinic subsequently extended this to $13$ lines \cite{DIM}, and Barakat and K\"uhne settled the case of up to $14$ lines over fields of arbitrary characteristic \cite{BK}. Most recently, Anghel and Stan proved the conjecture for arrangements of up to $16$ lines \cite{AS}. By contrast, to the best of our knowledge, much less is known for arrangements in higher-dimensional projective spaces. The main aim of this note is to take a step in this direction by studying arrangements of planes in $\mathbb{P}^3_{\mathbb{C}}$. Our main result that we would like to report is the following.

\begin{theoremA}
Let $\A,\mathcal{B}\subset\PP^{3}_{\C}$ be two arrangements of $d\le 9$ planes such that $L(\A)\cong L(\mathcal{B})$. If $\A$ is free, then $\mathcal{B}$ is free.
\end{theoremA}

In fact, we determine all intersection lattices of free arrangements of at most nine planes in $\PP^3_{\C}$.

\begin{theoremB}
Let $\A\subset\PP^{3}_{\C}$ be an essential and irreducible arrangement of $d\le9$ planes. Then $\A$ is free if and only if $L(\A)\cong L(\A_i)$ for one of the arrangements $\A_1,\ldots,\A_{19}$ in Table~\ref{tab:list}. In particular:
\begin{enumerate}
\item $d\in\{7,8,9\}$ and exponents $\exps(\A)\in\{(2,2,2),(2,2,3),(2,3,3),(2,2,4)\}$;
\item $\A$ is free if and only if $\A$ is divisionally free.
\end{enumerate}
\end{theoremB}

Our approach is based on two simple observations. First, an irreducible free arrangement cannot have $1$ among its exponents (Lemma~\ref{lem:irr}), which forces $d\ge 7$. Second, by Ziegler's restriction theorem, if $\A$ is free with exponents $d_1\le d_2\le d_3$, then each plane $H\in\A$ meets the remaining planes along at most $d-d_1$ distinct lines (Lemma~\ref{lem:restr}). Combining this with a simple count involving the numbers $t_k$ of lines contained in exactly $k$ planes, we obtain a short list of admissible weak combinatorics (Proposition~\ref{prop:weak}). It is worth emphasizing that the computer plays only a marginal role in our considerations: it is used once, in Proposition~\ref{prop:list}, to list the intersection lattices with the admissible weak combinatorics by inspecting the database of matroids constructed in \cite{MMIB}. All remaining arguments are by hand.

As a by-product of Theorem B, we show that no Cynk--Szemberg octic arrangement is free. These are the arrangements of eight planes that are used to construct double octic Calabi--Yau threefolds, see Corollary~\ref{cor:octic}.

In the last section, motivated by the recent work of Dimca and Pokora on Ziegler pairs of plane arrangements \cite{DP}, we discuss to what extent the minimal degree of a Jacobian relation of a plane arrangement is determined by its combinatorics. We present examples of weak Ziegler pairs, i.e., pairs of irreducible arrangements with the same numbers of multiple lines and multiple points, but with different minimal degrees of Jacobian relations, including pairs of octic arrangements leading to rigid double octics from the classification of Cynk and Kocel-Cynk \cite{CKC}, and we explain why such phenomena are hard to find among arrangements with isomorphic intersection lattices.

Throughout the note we work over the complex numbers.

\section{Preliminaries}

Let $\A=\{H_1,\ldots,H_d\}\subset\PP^3_{\C}$ be an arrangement of $d$ planes with $H_i=V(\ell_i)$, and let $S=\C[x,y,z,w]$. We say that $\A$ is \emph{essential} if $H_1\cap\ldots\cap H_d=\emptyset$, and \emph{irreducible} if there is no linear change of coordinates after which $Q(\A)$ becomes a product of two non-constant polynomials in disjoint sets of variables. If $\A$ is not essential, then $\A$ is a cone over a line arrangement. If $\A$ is essential and reducible, then $\A=\A'\times\A''$ is a product of two arrangements of smaller rank. Both properties can be read off from the intersection lattice $L(\A)$: in matroid language, $\A$ is irreducible if and only if its matroid is connected.

A line $L\subset\PP^3_{\C}$ is called a \emph{$k$-fold line} of $\A$ if exactly $k\ge 2$ planes of $\A$ contain $L$, and we denote by $t_k=t_k(\A)$ the number of $k$-fold lines. Lines with $k\ge3$ are called \emph{multiple}. For a plane $H\in\A$ the \emph{restriction}
\[
\A^H=\{H\cap H'\,:\,H'\in\A\setminus\{H\}\}
\]
is a line arrangement in $H\cong\PP^2_{\C}$, and we put $n_H=|\A^H|$. Similarly, for a line $L=H\cap H'$ we denote by $\A^L$ the set of points $L\cap H''$ with $H''\in\A$, $L\not\subset H''$.

We denote by $\chi(\A;t)$ the characteristic polynomial of the central arrangement in $\C^4$ defined by $Q(\A)$, and we set $\chi_0(\A;t)=\chi(\A;t)/(t-1)$. Writing
\[
\chi(\A;t)=t^4-b_1t^3+b_2t^2-b_3t+b_4,
\]
we have $b_1=d$ and $b_2=\sum_{k\ge2}(k-1)\,t_k$. Recall that $\chi(\A;t)$ is determined by $L(\A)$.

Every derivation $\theta\in\Der(S)$ splits as $\theta=f\theta_E+\theta_0$ with $\theta_E=x\partial_x+y\partial_y+z\partial_z+w\partial_w$ the Euler derivation and $\theta_0(Q(\A))=0$, so $D(\A)=S\cdot\theta_E\oplus D_0(\A)$, where 
$$D_0(\A)=\{\theta\in\Der(S):\theta(Q(\A))=0\}.$$ If $\A$ is free, then $D_0(\A)\cong S(-d_1)\oplus S(-d_2)\oplus S(-d_3)$, and we say that $\A$ is free with exponents $\exps(\A)=(d_1,d_2,d_3)$, where $d_1\le d_2\le d_3$ and $d_1+d_2+d_3=d-1$, see \cite{OT}.

\begin{theorem}[Terao's factorization {\cite{Te}}]\label{thm:fact}
If $\A$ is free with exponents $(d_1,d_2,d_3)$, then $\chi_0(\A;t)=(t-d_1)(t-d_2)(t-d_3)$.
\end{theorem}

Next we recall the notion of Ziegler restrictions. For $H\in\A$ and $K\in\A^H$ put $m^H(K)=\#\{H'\in\A\setminus\{H\}: H\cap H'=K\}$. Then $(\A^H,m^H)$ is a multiarrangement in $H\cong\PP^2_{\C}$ with $\sum_{K\in\A^H}m^H(K)=d-1$. If $\alpha_K$ is a linear form defining $K$ in the coordinate ring $S_H$ of $H$, we set
\[
D(\A^H,m^H)=\{\theta\in\Der(S_H)\,:\,\theta(\alpha_K)\in\alpha_K^{m^H(K)}\cdot S_H \text{ for all } K\in\A^H\}.
\]

\begin{theorem}[Ziegler {\cite{Z}}]\label{thm:ziegler}
If $\A$ is free with exponents $(d_1,d_2,d_3)$, then for every $H\in\A$ the multiarrangement $(\A^H,m^H)$ is free with exponents $(d_1,d_2,d_3)$.
\end{theorem}

Now we recall Abe's division theorem and the notion of divisional freeness.

\begin{theorem}[Abe {\cite{Abe}}]\label{thm:division}
Let $H\in\A$. If $\A^H$ is free and $\chi_0(\A^H;t)$ divides $\chi_0(\A;t)$, then $\A$ is free.
\end{theorem}

\begin{definition}[{\cite{Abe}}]
An essential arrangement $\A\subset\PP^3_{\C}$ is \emph{divisionally free} if there exist a plane $H\in\A$ and a line $L=H\cap H'$ with $H'\in\A$ such that
\[
\chi_0(\A^H;t)\ \text{divides}\ \chi_0(\A;t)\qquad\text{and}\qquad \chi_0(\A^L;t)=t-(|\A^L|-1)\ \text{divides}\ \chi_0(\A^H;t).
\]
\end{definition}

Since every arrangement of points on a line is free, Theorem~\ref{thm:division} applied twice shows that a divisionally free arrangement is free. Moreover, divisional freeness depends only on $L(\A)$. These two facts are due to Abe \cite{Abe}.

Finally, we need the following classical facts. If $\A=\A'\times\A''$ is a product, then $\A$ is free if and only if both $\A'$ and $\A''$ are free \cite[Proposition 4.28]{OT}. If $\A$ is a cone over a line arrangement $\mathcal{C}$, then $\A$ is free if and only if $\mathcal{C}$ is free. 
\section{Restrictions on the weak combinatorics of free arrangements}

We start with the following simple lemma.

\begin{lemma}\label{lem:irr}
Let $\A\subset\PP^{3}_{\C}$ be an essential and irreducible arrangement. Then
\[
D(\A)_1=\C\cdot\theta_E .
\]
In particular, if $\A$ is free with exponents $(d_1,d_2,d_3)$, then $d_1\ge 2$.
\end{lemma}

\begin{proof}
We write $\mathbf{x}=(x,y,z,w)^{T}$, and for a linear form $\ell$ we denote by $a\in\C^4$ its coefficient vector, so that $\ell=a^{T}\mathbf{x}$. Let $a_i$ be the coefficient vector of $\ell_i$. Every derivation of degree $1$ is of the form
\[
\theta_M=\sum_{j=1}^{4}(M\mathbf{x})_j\,\partial_j \qquad\text{for a unique } M\in{\rm Mat}_{4\times4}(\C),
\]
and $\theta_E=\theta_{\rm Id}$. For $\ell=a^{T}\mathbf{x}$ we have
\[
\theta_M(\ell)=a^{T}M\mathbf{x}=(M^{T}a)^{T}\mathbf{x},
\]
which is again a linear form. We proceed in three steps.

\emph{Step 1: $\theta_M\in D(\A)$ if and only if every $a_i$ is an eigenvector of $M^{T}$.} Since $\theta_M(\ell_i)$ is a linear form, it is divisible by $\ell_i$ if and only if $\theta_M(\ell_i)=\lambda_i\ell_i$ for some $\lambda_i\in\C$, that is, $M^{T}a_i=\lambda_i a_i$.

\emph{Step 2: If $\theta_M\in D(\A)$, then all $\lambda_i$ coincide.} Suppose that $\mu_1,\ldots,\mu_s$, with $s\ge2$, are the distinct values among $\lambda_1,\ldots,\lambda_d$. For $j=1,\ldots,s$ put
\[
\A_j=\{H_i\in\A\,:\,\lambda_i=\mu_j\},\qquad W_j={\rm span}\{a_i\,:\,H_i\in\A_j\}\subset\C^4.
\]
Each $W_j$ is contained in the $\mu_j$-eigenspace of $M^{T}$. Eigenspaces corresponding to distinct eigenvalues are linearly independent, so the sum $W_1+\ldots+W_s$ is direct. Since $\A$ is essential, the vectors $a_1,\ldots,a_d$ span $\C^4$, and therefore
\[
\C^4=W_1\oplus\ldots\oplus W_s .
\]
Now take as new coordinates a basis of linear forms adapted to this decomposition. Then every form defining a plane of $\A_j$ involves only the coordinates coming from $W_j$, so
\[
Q(\A)=Q(\A_1)\cdots Q(\A_s)
\]
with factors in pairwise disjoint sets of variables. Hence $\A$ is reducible, a contradiction.

\emph{Step 3: $\theta_M=\lambda\,\theta_E$.} By Step 2 there is $\lambda\in\C$ with $M^{T}a_i=\lambda a_i$ for all $i$. The vectors $a_i$ span $\C^4$, so $M^{T}=\lambda\,{\rm Id}$, i.e., $\theta_M=\lambda\,\theta_E$. This proves that $D(\A)_1=\C\cdot\theta_E$.

Finally, recall that $D(\A)=S\cdot\theta_E\oplus D_0(\A)$. If $\A$ were free with $d_1=1$, then $D_0(\A)$ would contain a non-zero derivation of degree $1$. Since the sum is direct, this derivation is not a multiple of $\theta_E$, which contradicts $D(\A)_1=\C\cdot\theta_E$. Hence $d_1\ge2$.
\end{proof}

\begin{corollary}\label{cor:exps}
Let $\A\subset\PP^3_{\C}$ be an essential and irreducible free arrangement of $d\le9$ planes. Then $d\in\{7,8,9\}$ and
\[
(d;\exps(\A))\in\{(7;2,2,2),\ (8;2,2,3),\ (9;2,3,3),\ (9;2,2,4)\}.
\]
\end{corollary}

\begin{proof}
By Lemma~\ref{lem:irr} we have $2\le d_1\le d_2\le d_3$, and $d_1+d_2+d_3=d-1\le8$. The claim follows.
\end{proof}

The next lemma is a key ingredient of our proof.

\begin{lemma}\label{lem:restr}
Let $\A\subset\PP^3_{\C}$ be a free arrangement of $d$ planes with exponents $d_1\le d_2\le d_3$. Then for every plane $H\in\A$ one has
\[
n_H\le d-d_1 .
\]
\end{lemma}

\begin{proof}
By Theorem~\ref{thm:ziegler}, the Ziegler restriction $(\A^H,m^H)$ is free with exponents $(d_1,d_2,d_3)$. Consider
\[
\theta=\Big(\prod_{K\in\A^H}\alpha_K^{m^H(K)-1}\Big)\cdot\theta_E\in\Der(S_H).
\]
Since $\theta(\alpha_K)$ is divisible by $\alpha_K^{m^H(K)}$ for every $K$, we have $0\ne\theta\in D(\A^H,m^H)$, and
\[
\deg\theta=1+\sum_{K\in\A^H}\big(m^H(K)-1\big)=1+(d-1)-n_H=d-n_H.
\]
A non-zero homogeneous element of a free graded module has degree at least the smallest degree of a basis element, hence $d_1\le d-n_H$.
\end{proof}

\begin{corollary}\label{cor:cover}
Let $\A\subset\PP^3_{\C}$ be an essential and irreducible free arrangement. Then every plane $H\in\A$ contains a multiple line of $\A$. In particular,
\begin{equation}\label{eq:cover}
\sum_{k\ge3}k\,t_k\ \ge\ d.
\end{equation}
\end{corollary}

\begin{proof}
Every $k$-fold line $L\subset H$ accounts for $k-1$ planes $H'\ne H$ but only for one element of $\A^H$, so
\[
n_H=(d-1)-\sum_{L\subset H}(k_L-2),
\]
where the sum runs over the multiple lines $L\subset H$ and $k_L$ denotes the multiplicity of $L$. By Lemmas~\ref{lem:irr} and~\ref{lem:restr}, $n_H\le d-2$, so the sum is positive. Since every plane contains a multiple line, counting the pairs (plane, multiple line in it) gives \eqref{eq:cover}.
\end{proof}

We also need the following elementary count.

\begin{lemma}\label{lem:count}
Let $\A\subset\PP^{3}_{\C}$ be an arrangement of $d$ planes. Then
\[
\binom{d}{2}=\sum_{k\ge2}\binom{k}{2}t_k,\qquad\text{hence}\qquad\sum_{k\ge3}\binom{k-1}{2}\,t_k=\binom{d}{2}-b_2.
\]
\end{lemma}

\begin{proof}
Every pair of distinct planes meets along exactly one line, which gives the first equality. Subtracting $b_2=\sum_{k\ge2}(k-1)t_k$ and using $\binom{k}{2}-(k-1)=\binom{k-1}{2}$ gives the second.
\end{proof}

Now we are ready to determine the admissible weak combinatorics.

\begin{proposition}\label{prop:weak}
Let $\A\subset\PP^{3}_{\C}$ be an essential and irreducible free arrangement of $d\le9$ planes. Then $t_k=0$ for $k\ge6$, and $(d;\exps(\A);t_3,t_4,t_5)$ is one of the following:
\begin{center}
\begin{tabular}{cccl}
\toprule
$d$ & $\exps(\A)$ & $\binom{d}{2}-b_2$ & $(t_3,t_4,t_5)$\\
\midrule
$7$ & $(2,2,2)$ & $3$ & $(3,0,0)$\\
$8$ & $(2,2,3)$ & $5$ & $(5,0,0)$, $(2,1,0)$\\
$9$ & $(2,3,3)$ & $7$ & $(7,0,0)$, $(4,1,0)$, $(1,2,0)$\\
$9$ & $(2,2,4)$ & $8$ & $(8,0,0)$, $(5,1,0)$, $(2,2,0)$, $(2,0,1)$\\
\bottomrule
\end{tabular}
\end{center}
\end{proposition}

\begin{proof}
By Corollary~\ref{cor:exps} and Theorem~\ref{thm:fact}, $b_2=(d-1)+d_1d_2+d_1d_3+d_2d_3$, which gives the third column. By Lemma~\ref{lem:count},
\[
t_3+3t_4+6t_5+10t_6+\ldots=\binom{d}{2}-b_2\le 8,
\]
so $t_k=0$ for $k\ge6$. Solving this equation in non-negative integers gives, besides the solutions listed above, exactly two further cases: $(t_3,t_4,t_5)=(0,1,0)$ for $d=7$, and $(t_3,t_4,t_5)=(1,0,1)$ for $d=9$ and $\exps(\A)=(2,3,3)$. In these two cases $\sum_k k\,t_k$ equals $4<7$ and $8<9$, respectively, which contradicts Corollary~\ref{cor:cover}.
\end{proof}

\begin{remark}
In the excluded case $(9;(2,3,3);1,0,1)$ the covering argument is completely transparent, namely the $5$-fold line and the triple line lie on at most $5+3=8$ planes, so some plane $H$ contains no multiple line, and then $n_H=8=d-1$. By Lemma~\ref{lem:restr} this would force $d_1=1$, contradicting Lemma~\ref{lem:irr}.
\end{remark}

\section{The intersection lattices}

The following is the only place where we use computer.

\begin{proposition}\label{prop:list}
Let $\A\subset\PP^{3}_{\C}$ be an essential and irreducible arrangement of $d\le9$ planes such that $\chi_0(\A;t)$ and $(t_3,t_4,t_5)$ are as in one of the rows of Proposition~\ref{prop:weak}, and every plane of $\A$ contains a multiple line. Then $L(\A)\cong L(\A_i)$ for exactly one of the arrangements $\A_1,\ldots,\A_{19}$ of Table~\ref{tab:list}.
\end{proposition}

\begin{proof}
The intersection lattice of $\A$ is the lattice of flats of a simple, connected matroid of rank $4$ on $d$ elements. We inspect the complete list of matroids of rank $4$ on at most $9$ elements due to Matsumoto, Moriyama, Imai and Bremner \cite{MMIB}, which is also available in \textsc{SageMath}. Among the simple and connected ones, exactly $2$, $4$ and $16$ matroids on $7$, $8$ and $9$ elements have $\chi_0$ as in Proposition~\ref{prop:weak}. Two of the sixteen matroids on nine elements have weak combinatorics $(t_3,t_4,t_5)=(1,0,1)$ and are excluded by Corollary~\ref{cor:cover}. One more matroid on nine elements, with $\exps=(2,2,4)$ and $t_3=8$, contains the Fano plane as a restriction, so it is not realizable over $\C$. The remaining $19$ matroids are realized by the arrangements $\A_1,\ldots,\A_{19}$. The search takes a couple of minutes on a standard laptop.
\end{proof}

\begin{table}[ht]
\centering
\small
\caption{The intersection lattices of essential and irreducible free arrangements of $d\le9$ planes in $\PP^3_{\C}$. The symbol $\dagger$ marks graphic arrangements, i.e., arrangements associated with chordal graphs on five vertices.}\label{tab:list}
\begin{tabular}{ccccl}
\toprule
 & $d$ & $\exps(\A)$ & $(t_3,t_4,t_5)$ & defining polynomial\\
\midrule
$\mathcal{A}_{1}$ & $7$ & $(2,2,2)$ & $(3,0,0)$ & $Q=xy(x+y)z(x+z)w(y+w)$\rlap{$^\dagger$} \\
$\mathcal{A}_{2}$ & $7$ & $(2,2,2)$ & $(3,0,0)$ & $Q=xy(x+y)z(x+z)w(x+w)$\rlap{$^\dagger$} \\
$\mathcal{A}_{3}$ & $8$ & $(2,2,3)$ & $(5,0,0)$ & $Q=xy(x+y)z(z-x)(y+z)w(x+w)$\rlap{$^\dagger$} \\
$\mathcal{A}_{4}$ & $8$ & $(2,2,3)$ & $(2,1,0)$ & $Q=xy(x+y)(y-x)z(x+z)w(z+w)$ \\
$\mathcal{A}_{5}$ & $8$ & $(2,2,3)$ & $(2,1,0)$ & $Q=xy(x+y)(y-x)z(x+z)w(y+w)$ \\
$\mathcal{A}_{6}$ & $8$ & $(2,2,3)$ & $(2,1,0)$ & $Q=xy(x+y)(y-x)z(x+z)w(x+w)$ \\
$\mathcal{A}_{7}$ & $9$ & $(2,3,3)$ & $(7,0,0)$ & $Q=xy(x+y)z(z-x)(y+z)w(w-x)(y+w)$\rlap{$^\dagger$} \\
$\mathcal{A}_{8}$ & $9$ & $(2,3,3)$ & $(1,2,0)$ & $Q=xy(x+y)(y-x)z(x+z)w(z-w)(z+w)$ \\
$\mathcal{A}_{9}$ & $9$ & $(2,3,3)$ & $(4,1,0)$ & $Q=xy(x+y)(y-x)z(x-z)w(x+w)(z+w)$ \\
$\mathcal{A}_{10}$ & $9$ & $(2,3,3)$ & $(7,0,0)$ & $Q=xy(x+y)z(x+z)(y-z)(x+y+z)w(y+w)$ \\
$\mathcal{A}_{11}$ & $9$ & $(2,3,3)$ & $(7,0,0)$ & $Q=xy(x+y)z(x+z)(y-z)(x+y+z)w(x+w)$ \\
$\mathcal{A}_{12}$ & $9$ & $(2,3,3)$ & $(4,1,0)$ & $Q=xy(x+y)(y-x)z(z-x)(y+z)w(z+w)$ \\
$\mathcal{A}_{13}$ & $9$ & $(2,3,3)$ & $(4,1,0)$ & $Q=xy(x+y)(y-x)z(z-x)(y+z)w(y+w-x)$ \\
$\mathcal{A}_{14}$ & $9$ & $(2,3,3)$ & $(4,1,0)$ & $Q=xy(x+y)(y-x)z(z-x)(y+z)w(x+w)$ \\
$\mathcal{A}_{15}$ & $9$ & $(2,3,3)$ & $(1,2,0)$ & $Q=xy(x+y)(y-x)z(x-z)(x+z)w(y+w)$ \\
$\mathcal{A}_{16}$ & $9$ & $(2,3,3)$ & $(1,2,0)$ & $Q=xy(x+y)(y-x)z(x-z)(x+z)w(x+w)$ \\
$\mathcal{A}_{17}$ & $9$ & $(2,2,4)$ & $(2,0,1)$ & $Q=xy(x+y)(y-x)(2x+y)z(x+z)w(z+w)$ \\
$\mathcal{A}_{18}$ & $9$ & $(2,2,4)$ & $(2,0,1)$ & $Q=xy(x+y)(y-x)(2x+y)z(x+z)w(y+w)$ \\
$\mathcal{A}_{19}$ & $9$ & $(2,2,4)$ & $(2,0,1)$ & $Q=xy(x+y)(y-x)(2x+y)z(x+z)w(x+w)$ \\
\bottomrule
\end{tabular}
\end{table}

\begin{lemma}\label{lem:DF}
Each of the arrangements $\A_1,\ldots,\A_{19}$ is divisionally free, and hence free.
\end{lemma}

\begin{proof}
For each $\A_i$, Table~\ref{tab:DF} gives a plane $H\in\A_i$, the number $n_H$, the exponents of the line arrangement $\A_i^H$ (that is, the roots of $\chi_0(\A_i^H;t)$), a second plane $H'$, and the number of points on the line $L=H\cap H'$. In every case the exponents of $\A_i^H$ are contained in $\exps(\A_i)$, and $|\A_i^L|-1$ is one of the exponents of $\A_i^H$, which is exactly the definition of divisional freeness. The verification is a direct inspection of the multiple points of $\A_i^H$.

For instance, take $\A_1$ with $Q=xy(x+y)z(x+z)w(y+w)$ and $H=\{x+y=0\}$. On $H$ the planes $x=0$, $y=0$, $x+y=0$ all cut out the line $\{x=0\}$, so $\A_1^H$ consists of the five lines $x=0$, $z=0$, $x+z=0$, $w=0$, $w-x=0$. These lines have two triple points, $\{x=z=0\}$ and $\{x=w=0\}$, and four double points, so
\[
\chi_0(\A_1^H;t)=t^2-4t+4=(t-2)^2,
\]
which divides $\chi_0(\A_1;t)=(t-2)^3$. The line $L=H\cap\{z=0\}$ meets $\A_1^H$ in the three points $\{x=z=0\}$, $\{z=w=0\}$ and $\{z=0,\ w=x\}$, so $\chi_0(\A_1^L;t)=t-2$ divides $(t-2)^2$. The remaining cases are analogous.
\end{proof}

\begin{table}[ht]
\centering
\small
\caption{Divisional freeness: $L=H\cap H'$.}\label{tab:DF}
\begin{tabular}{cccccc}
\toprule
 & $H$ & $n_H$ & $\exps(\A^H)$ & $H'$ & $|\A^L|$\\
\midrule
$\mathcal{A}_{1}$ & $x+y=0$ & $5$ & $(2,2)$ & $z=0$ & $3$ \\
$\mathcal{A}_{2}$ & $y=0$ & $5$ & $(2,2)$ & $z=0$ & $3$ \\
$\mathcal{A}_{3}$ & $y=0$ & $5$ & $(2,2)$ & $z=0$ & $3$ \\
$\mathcal{A}_{4}$ & $y=0$ & $5$ & $(2,2)$ & $x=0$ & $3$ \\
$\mathcal{A}_{5}$ & $x+y=0$ & $5$ & $(2,2)$ & $z=0$ & $3$ \\
$\mathcal{A}_{6}$ & $y=0$ & $5$ & $(2,2)$ & $z=0$ & $3$ \\
$\mathcal{A}_{7}$ & $z=0$ & $6$ & $(2,3)$ & $x=0$ & $3$ \\
$\mathcal{A}_{8}$ & $y=0$ & $6$ & $(2,3)$ & $x=0$ & $4$ \\
$\mathcal{A}_{9}$ & $y=0$ & $6$ & $(2,3)$ & $x=0$ & $3$ \\
$\mathcal{A}_{10}$ & $x=0$ & $6$ & $(2,3)$ & $z=0$ & $3$ \\
$\mathcal{A}_{11}$ & $y-z=0$ & $6$ & $(2,3)$ & $y=0$ & $3$ \\
$\mathcal{A}_{12}$ & $y-x=0$ & $6$ & $(2,3)$ & $x=0$ & $3$ \\
$\mathcal{A}_{13}$ & $z=0$ & $6$ & $(2,3)$ & $x=0$ & $3$ \\
$\mathcal{A}_{14}$ & $y-x=0$ & $6$ & $(2,3)$ & $z=0$ & $3$ \\
$\mathcal{A}_{15}$ & $x+y=0$ & $6$ & $(2,3)$ & $z=0$ & $3$ \\
$\mathcal{A}_{16}$ & $y=0$ & $6$ & $(2,3)$ & $z=0$ & $3$ \\
$\mathcal{A}_{17}$ & $y=0$ & $5$ & $(2,2)$ & $x=0$ & $3$ \\
$\mathcal{A}_{18}$ & $x+y=0$ & $5$ & $(2,2)$ & $z=0$ & $3$ \\
$\mathcal{A}_{19}$ & $y=0$ & $5$ & $(2,2)$ & $z=0$ & $3$ \\
\bottomrule
\end{tabular}
\end{table}

\section{Proofs of the main results}

\begin{proof}[Proof of Theorem B]
Let $\A$ be essential, irreducible and free with $d\le9$. By Corollary~\ref{cor:exps}, Theorem~\ref{thm:fact}, Corollary~\ref{cor:cover} and Proposition~\ref{prop:weak}, $\A$ satisfies the assumptions of Proposition~\ref{prop:list}, so $L(\A)\cong L(\A_i)$ for some $i$. Conversely, if $L(\A)\cong L(\A_i)$, then $\A$ is divisionally free by Lemma~\ref{lem:DF}, since divisional freeness depends only on the intersection lattice, and hence $\A$ is free by Theorem~\ref{thm:division}. Statements (1) and (2) follow from Table~\ref{tab:list} and Lemma~\ref{lem:DF}.
\end{proof}

\begin{proof}[Proof of Theorem A]
Let $\A,\mathcal{B}\subset\PP^3_{\C}$ be arrangements of $d\le9$ planes with $L(\A)\cong L(\mathcal{B})$, and assume that $\A$ is free. Being essential and being irreducible are properties of the intersection lattice, so $\A$ and $\mathcal{B}$ fall into the same of the following three cases.

\emph{Case 1: $\A$ is not essential.} Then $\A$ and $\mathcal{B}$ are cones over line arrangements with $d\le9$ lines and isomorphic intersection lattices, and the claim follows from \cite{FV}.

\emph{Case 2: $\A$ is essential and reducible.} Then $\A=\A'\times\A''$ and $\mathcal{B}=\mathcal{B}'\times\mathcal{B}''$ with $L(\A')\cong L(\mathcal{B}')$ and $L(\A'')\cong L(\mathcal{B}'')$, where the factors have rank at most $3$ and at most $8$ hyperplanes. Since $\A$ is free, so are $\A'$ and $\A''$. Arrangements of rank at most $2$ are always free, and in rank $3$ we use \cite{FV} again. Hence $\mathcal{B}'$ and $\mathcal{B}''$ are free, and so is $\mathcal{B}$.

\emph{Case 3: $\A$ is essential and irreducible.} By Theorem B, $L(\mathcal{B})\cong L(\A)\cong L(\A_i)$ for some $i$, and again by Theorem B the arrangement $\mathcal{B}$ is free.
\end{proof}

We should point out here that the factorization of the characteristic polynomial alone does not imply freeness, even in our small range. The following example shows the two lattices excluded in Proposition~\ref{prop:list}.

\begin{example}
Consider the arrangements
\begin{align*}
\mathcal{N}_1:&\quad Q=xy(x+y)(y-x)(2x+y)\,z\,(y+z-2x)\,w\,(z+w),\\
\mathcal{N}_2:&\quad Q=xy(x+y)(y-x)(2x+y)\,z\,(x+z)\,w\,(y+w-2x).
\end{align*}
In both cases $d=9$, $\chi_0(t)=(t-2)(t-3)^2$, and the weak combinatorics is $(t_3,t_4,t_5)=(1,0,1)$: the planes $x,y,x+y,y-x,2x+y$ pass through the $5$-fold line $\{x=y=0\}$, and the triple line is $\{z=w=0\}$ in $\mathcal{N}_1$ and $\{x=z=0\}$ in $\mathcal{N}_2$. The plane $y+z-2x=0$ in $\mathcal{N}_1$, and the planes $w=0$ and $y+w-2x=0$ in $\mathcal{N}_2$, contain no multiple line. Hence neither arrangement is free, and the same holds for every arrangement with one of these intersection lattices. One can check with \textsc{Singular} \cite{Sing} that $D_0(\mathcal{N}_i)$ is minimally generated in degrees $2,2,5,5$. A direct reason is that the derivations
\[
w(z+w)\,\partial_w,\qquad z(y+z-2x)\,\partial_z+w(y-w-2x)\,\partial_w
\]
both belong to $D(\mathcal{N}_1)$ and are linearly independent modulo $S_1\cdot\theta_E$, so $\dim D_0(\mathcal{N}_1)_2\ge2$, whereas every free arrangement with exponents $(2,3,3)$ satisfies $\dim D_0(\A)_2=1$.
\end{example}

\begin{remark}
The same argument works for line arrangements. If $\mathcal{C}\subset\PP^2_{\C}$ is a free arrangement of $d$ lines with exponents $d_1\le d_2$, then $n_H\le d-d_1$ for every line $H\in\mathcal{C}$. For instance, the arrangement of five lines through one point and two further general lines has $\chi_0(t)=(t-3)^2$, but each of the two general lines meets the other lines in $6>7-3$ points, so it is never free. Note that its cone is a localization of $\mathcal{N}_1$ and of $\mathcal{N}_2$.
\end{remark}

Theorem B has a nice consequence for arrangements of eight planes that are important in the construction of Calabi--Yau threefolds. Recall that an arrangement of eight planes in $\PP^3_{\C}$ is called a \emph{Cynk--Szemberg octic arrangement} if it has only double and triple lines, and each point of $\PP^3_{\C}$ is contained in at most five planes of the arrangement \cite[Definition 3.3]{JP}. The double cover of $\PP^3_{\C}$ branched along such an arrangement has a crepant resolution which is a Calabi--Yau threefold \cite{CS, Me}, see Section~\ref{sec:weak} for more details.

\begin{corollary}\label{cor:octic}
A Cynk--Szemberg octic arrangement is never free.
\end{corollary}
\begin{proof}
Let $\A$ be a free arrangement of eight planes. If $\A$ is not essential, then all eight planes pass through one point. If $\A$ is essential and reducible, then after a linear change of coordinates $Q(\A)=Q'(u_1,\ldots,u_a)\cdot Q''(u_{a+1},\ldots,u_4)$ with $1\le a\le 3$. If $a=1$, then $Q'=u_1$ up to a constant, and the seven planes defined by the factors of $Q''$ pass through the point $(1:0:0:0)$. The case $a=3$ is symmetric. If $a=2$, then the planes defined by the factors of $Q'$ contain the line $\{u_1=u_2=0\}$, and the planes defined by the factors of $Q''$ contain the line $\{u_3=u_4=0\}$, so one of these lines is contained in at least four planes of $\A$. Finally, if $\A$ is essential and irreducible, then by Theorem B we have $L(\A)\cong L(\A_i)$ for some $i\in\{3,4,5,6\}$. The arrangements $\A_4,\A_5,\A_6$ have a quadruple line, and in $\A_3$ the six planes $x,y,x+y,z,z-x,y+z$ pass through the point $(0:0:0:1)$. In all cases $\A$ is not a Cynk--Szemberg octic arrangement.
\end{proof}

This result agrees with the observations made in \cite{JP}, where it is shown that the rigid Cynk--Szemberg octic arrangements are not free, and that general members of the $63$ one-parameter families of such arrangements from \cite{Me} are neither free nor nearly free. Moreover, the free arrangement $\mathcal{D}:xyzw(x+y)(y+z)(x-w)(x+y-w)=0$ found in \cite[Example 6.1]{JP} as a degeneration of a one-parameter family of octic arrangements has the same intersection lattice as $\A_3$, and it has a sixfold point.

It would be very interesting to extend Theorems A and B to arrangements of $d\ge10$ planes. In this range a complete list of matroids of rank $4$ is not available, so one would need a construction of the admissible intersection lattices that takes into account Lemma~\ref{lem:restr} and Lemma~\ref{lem:count} from the very beginning.

\section{Weak Ziegler pairs of plane arrangements}\label{sec:weak}

In this section we discuss a question which is closely related to Terao's conjecture, namely to what extent the algebraic properties of the module $D_0(\A)$ are determined by the combinatorics of $\A$. For line arrangements the first striking phenomenon of this kind was observed by Ziegler \cite{Z2}. He found two arrangements of $9$ lines with isomorphic intersection lattices, both consisting of the six sides and the three main diagonals of a hexagon, such that the six vertices of the hexagon lie on a conic in the first case, but not in the second one. The minimal degrees of the Jacobian relations of these two arrangements are equal to $5$ and $6$, respectively, see also \cite{DS}. Recently, Dimca and Pokora \cite{DP} started a systematic study of this phenomenon for arrangements of planes in $\PP^3_{\C}$. Our aim here is to present a weaker, purely numerical version of this phenomenon, which can be observed already for arrangements of $8$ planes, and in particular for octic arrangements defining rigid Calabi--Yau threefolds.

Let $\A\subset\PP^3_{\C}$ be an arrangement of $d$ planes and let $f=Q(\A)\in S=\C[x,y,z,w]$ be its defining polynomial. We denote by $J_f$ the Jacobian ideal of $f$, i.e., the ideal spanned by the partial derivatives $f_x,f_y,f_z,f_w$, and by $M(f)=S/J_f$ the Jacobian algebra of $f$. Recall that
\[
D_0(\A)=\{\theta=a\partial_x+b\partial_y+c\partial_z+e\partial_w\,:\,\theta(f)=a f_x+b f_y+c f_z+e f_w=0\},
\]
so $D_0(\A)$ is exactly the module of Jacobian relations of $f$. Following \cite{DP}, we write the minimal free resolution of $D_0(\A)$ in the form
\[
0\longrightarrow\bigoplus_{k=1}^{r}S(-b_k)\longrightarrow\bigoplus_{j=1}^{q}S(-c_j)\longrightarrow\bigoplus_{i=1}^{p}S(-d_i)\longrightarrow D_0(\A)\longrightarrow0,
\]
and we call the increasing sequences $\mathbf{d}=(d_1,\ldots,d_p)$, $\mathbf{c}=(c_1,\ldots,c_q)$ and $\mathbf{b}=(b_1,\ldots,b_r)$ the \emph{graded Betti numbers} of the Jacobian algebra $M(f)$. As in \cite{DP}, we write $k_m$ if the number $k$ appears $m$ times, for instance $(3,4_2,5)=(3,4,4,5)$. The sequence $\mathbf{d}$ gives the degrees of a minimal system of generators of $D_0(\A)$.

\begin{definition}
The \emph{minimal degree of a Jacobian relation} of $\A$ is the integer
\[
\operatorname{mdr}(\A)=d_1=\min\{r\,:\,D_0(\A)_r\neq0\}.
\]
\end{definition}

If $\A$ is free with exponents $(d_1,d_2,d_3)$, then $\mathbf{d}=(d_1,d_2,d_3)$ and $\operatorname{mdr}(\A)=d_1$. Now we recall the main definition of \cite{DP}.

\begin{definition}[{\cite[Definition 1.1]{DP}}]\label{def:ZP}
Two arrangements $\A:f=0$ and $\A':f'=0$ in $\PP^3_{\C}$ form a \emph{Ziegler pair} if $L(\A)\cong L(\A')$ and the graded Betti numbers of $M(f)$ and $M(f')$ are different. We say that a Ziegler pair satisfies the condition (MDR) if $\operatorname{mdr}(\A)=\operatorname{mdr}(\A')$.
\end{definition}

Dimca and Pokora observed that Ziegler pairs of plane arrangements can be easily obtained from Ziegler pairs of line arrangements. If $\mathcal{B}:g=0$ and $\mathcal{B}':g'=0$ is a Ziegler pair of line arrangements, where $g,g'\in\C[x,y,z]$, then the arrangements $wg=0$ and $wg'=0$ form a Ziegler pair of plane arrangements \cite[Theorem 1.3, Example 3.3]{DP}, and for Ziegler's hexagon even the Hilbert polynomials of the Jacobian algebras are different \cite[Example 3.1]{DP}. Let us point out that the surface $wg=0$, which is called the cone over the curve $g=0$ in \cite{DP}, contains the plane $w=0$, so it is an essential and reducible arrangement in our terminology. By \cite[Theorem 1.6]{DP} one has $\operatorname{mdr}=1$ for such arrangements, hence these Ziegler pairs satisfy the condition (MDR) for trivial reasons. In fact, this is a general phenomenon.

\begin{remark}\label{rem:mdr12}
For essential arrangements small values of $\operatorname{mdr}$ have a very simple meaning. A constant derivation $\partial_v$ belongs to $D(\A)$ only if the point $v$ lies on all planes of $\A$, hence $\operatorname{mdr}(\A)\ge1$ for every essential arrangement. If $\A$ is essential and reducible, then after a linear change of coordinates $Q(\A)=Q'(u_1,\ldots,u_a)\cdot Q''(u_{a+1},\ldots,u_4)$, and the derivation $u_{a+1}\partial_{u_{a+1}}+\ldots+u_4\partial_{u_4}$ belongs to $D(\A)$, but it is not a multiple of $\theta_E$. Thus $\operatorname{mdr}(\A)=1$. On the other hand, if $\A$ is essential and irreducible, then $\operatorname{mdr}(\A)\ge2$ by Lemma~\ref{lem:irr}. Summing up, an essential arrangement $\A$ satisfies $\operatorname{mdr}(\A)=1$ if and only if $\A$ is reducible. This is the reason why below we restrict our attention to essential and irreducible arrangements.
\end{remark}

Dimca and Pokora also constructed an essential and irreducible Ziegler pair of $11$ planes \cite[Example 4.1]{DP}, which we discuss in Example~\ref{ex:DP} below. This pair also satisfies the condition (MDR). As far as we know, no Ziegler pair of essential and irreducible plane arrangements that does not satisfy the condition (MDR) is known. This motivates the following weaker notion.

We say that a point $P\in\PP^3_{\C}$ is an \emph{$m$-fold point} of $\A$ if exactly $m\ge3$ planes of $\A$ pass through $P$ and the intersection of these planes is the point $P$. We denote by $p_m=p_m(\A)$ the number of $m$-fold points of $\A$.

\begin{definition}
The \emph{weak combinatorics} of an arrangement $\A\subset\PP^3_{\C}$ of $d$ planes is the vector
\[
W(\A)=(d;\,t_2,t_3,\ldots;\,p_3,p_4,\ldots).
\]
\end{definition}

Obviously, $W(\A)$ is determined by $L(\A)$, but two arrangements with different intersection lattices might have the same weak combinatorics.

\begin{definition}
Let $\A,\A'\subset\PP^3_{\C}$ be two essential and irreducible arrangements of $d$ planes. We say that $(\A,\A')$ is a \emph{weak Ziegler pair} if $W(\A)=W(\A')$ and $\operatorname{mdr}(\A)\neq\operatorname{mdr}(\A')$.
\end{definition}

Every Ziegler pair of essential and irreducible arrangements which does not satisfy the condition (MDR) is a weak Ziegler pair. In the definition of a weak Ziegler pair we ask for different values of $\operatorname{mdr}$, and not only for different graded Betti numbers, since for arrangements with different intersection lattices a difference in the graded Betti numbers is not surprising at all. A weak Ziegler pair shows that $\operatorname{mdr}$ is not determined by the numbers of multiple lines and multiple points.

Using the database \cite{MMIB} once again, for every intersection lattice of an essential and irreducible arrangement of $d\le9$ planes we computed $\operatorname{mdr}$ of a randomly chosen realization, and then we grouped the intersection lattices according to their weak combinatorics. Let us point out that one has to be a bit careful here, since some intersection lattices are realizable over $\C$, but not over $\mathbb{Q}$. For instance, this happens for the rigid octic arrangements $\mathfrak{A}$ and $\mathfrak{B}$ from \cite{CKC}, which are defined over $\mathbb{Q}(\sqrt{-3})$, see also \cite[Section 5]{JP}. For $d\le8$ we determined the realization spaces of all such lattices directly. For $d\le7$ we did not find any weak Ziegler pair. For $d=8$ there are exactly three weak combinatorics that are shared by intersection lattices with different values of $\operatorname{mdr}$, and for $d=9$ we found $27$ such weak combinatorics. For all three cases with $d=8$, and for $22$ cases with $d=9$, we found explicit realizations defined over $\mathbb{Q}$ and we verified all claims below using \textsc{Singular} \cite{Sing}.

Before we present our examples, let us recall a classical motivation for studying arrangements of eight planes. Following \cite[Definition 3.3]{JP}, see also \cite{CS, Me}, an arrangement $\A\subset\PP^3_{\C}$ of eight planes is called a \emph{Cynk--Szemberg octic arrangement} if $\A$ has only double and triple lines, and each point of $\PP^3_{\C}$ is contained in at most five planes of $\A$. For such an arrangement the double cover of $\PP^3_{\C}$ branched along $\A$ has a crepant resolution of singularities $\widetilde{X}$ which is a Calabi--Yau threefold, the so-called \emph{double octic} \cite{CS, Me}. Cynk and Kocel-Cynk classified all double octics with $h^{1,2}(\widetilde{X})\le1$ defined by arrangements of eight planes \cite{CKC}. In particular, there are exactly $14$ arrangements leading to rigid double octics, i.e., with $h^{1,2}(\widetilde{X})=0$, and their modules of derivations were studied by Janasz and Pokora \cite{JP}. Following \cite[Definition 2.6]{JP}, which generalizes the notion introduced in \cite{ADP} for plane curves, the \emph{type} of an essential arrangement $\A$ of $d$ planes with $\mathbf{d}=(d_1,\ldots,d_p)$ is the number
\[
t(\A)=d_1+d_2+d_3-d+1.
\]
It is easy to see that $\A$ is free if and only if $t(\A)=0$, see \cite[Example 2.7]{JP}, so the type measures how far $\A$ is from being free.

Observe that all weak Ziegler pairs of eight planes that we found consist of Cynk--Szemberg octic arrangements, since these arrangements have at most one triple line and no point lying on five planes. Hence each such pair gives a pair of double octic Calabi--Yau threefolds. Let us start with the most interesting case, namely the case of rigid double octics. Below, $\mathcal{O}_n$ denotes the arrangement with number $n$ in the classification \cite{CKC}, and we use the equations from \cite[Section 5]{JP}.

\begin{example}\label{ex:octic}
Consider the following three arrangements of eight planes:
\begin{align*}
\mathcal{O}_{239}:&\quad f=xyzw(x+y+z)(x+y+w)(x+z+w)(y+z+w),\\
\mathcal{O}_{240}:&\quad f=xyzw(x+y+z)(x+y-z+w)(x-y+z+w)(x-y-z-w),\\
\mathcal{O}_{241}:&\quad f=xyzw(x+y+z+w)(x+y-z-w)(y-z+w)(x+z-w).
\end{align*}
These arrangements have only double lines, and each of them has exactly $16$ triple points and $10$ quadruple points, so
\[
W(\mathcal{O}_{239})=W(\mathcal{O}_{240})=W(\mathcal{O}_{241})=(8;\,28;\,16,10).
\]
The intersection lattices of these three arrangements are pairwise non-isomorphic, and the corresponding double octics are rigid Calabi--Yau threefolds \cite{CKC}. On the other hand, the graded Betti numbers of their Jacobian algebras are as follows, see also \cite[Section 5]{JP}:
\begin{align*}
\mathcal{O}_{239},\ \mathcal{O}_{240}:&\quad \mathbf{d}=(4_{10}),\quad \mathbf{c}=(5_9),\quad \mathbf{b}=(6_2),\\
\mathcal{O}_{241}:&\quad \mathbf{d}=(3,4_6),\quad \mathbf{c}=(5_3,6_2),\quad \mathbf{b}=(7).
\end{align*}
In particular, $\operatorname{mdr}(\mathcal{O}_{241})=3$ and $\operatorname{mdr}(\mathcal{O}_{239})=\operatorname{mdr}(\mathcal{O}_{240})=4$, so $(\mathcal{O}_{241},\mathcal{O}_{239})$ and $(\mathcal{O}_{241},\mathcal{O}_{240})$ are weak Ziegler pairs of rigid octic arrangements. In terms of the type, we have $t(\mathcal{O}_{241})=4$ and $t(\mathcal{O}_{239})=t(\mathcal{O}_{240})=5$.

Moreover, there is exactly one more intersection lattice of eight planes with the weak combinatorics $(8;\,28;\,16,10)$. It is the intersection lattice of the one-parameter family
\begin{multline*}
\mathcal{O}_{242}(A:B):\quad f=xyzw(x+y+z)(x+z-w)\\
\times\big(Ax+(A+B)y-Bz+Bw\big)\big((A+B)x+(A+B)y+Bw\big)
\end{multline*}
from \cite{Me, CKC}, so the corresponding double octics have $h^{1,2}=1$. For instance, for $(A:B)=(1:1)$ we get
\[
\mathcal{O}_{242}(1:1):\quad f=xyzw(x+y+z)(x+z-w)(x+2y-z+w)(2x+2y+w).
\]
For all members of this family that we checked, the graded Betti numbers of the Jacobian algebra are the same as for $\mathcal{O}_{241}$. Summing up, the four intersection lattices with the weak combinatorics $(8;\,28;\,16,10)$ split into two groups of two lattices according to the value of $\operatorname{mdr}$, and within each group the graded Betti numbers of the arrangements above coincide.
\end{example}

\begin{example}\label{ex:weak8}
Now we present the remaining two weak Ziegler pairs of eight planes:
\begin{align*}
\A_{1}:&\quad f=xyz(x+y+z)\,w\,(y-2x+2w)(z-2x-2w)(z-y+w),\\
\A'_{1}:&\quad f=xyz(x+y+z)\,w\,(x-y+2w)(x+2z-2w)(z-y+w),\\[4pt]
\A_{2}:&\quad f=xy(x+y)zw\,(2w-3x-6z)(3y+6z+2w)(x+2y+2z+2w),\\
\A'_{2}:&\quad f=xy(x+y)zw\,(2w-x-z)(2y+z+4w)(x+y+z+w).
\end{align*}
All these arrangements are essential, irreducible, and they are Cynk--Szemberg octic arrangements. We have
\[
W(\A_1)=W(\A'_1)=(8;\,28;\,24,8),\qquad W(\A_2)=W(\A'_2)=(8;\,25,1;\,16,11).
\]
In Table~\ref{tab:weak8} we collect the graded Betti numbers, the types, and the dimensions of the realization spaces $\mathcal{R}$ of the corresponding intersection lattices for all weak Ziegler pairs of eight planes. In particular, the intersection lattices of $\A_1,\A'_1,\A_2,\A'_2$ have positive-dimensional realization spaces, so these arrangements are not rigid. Comparing the incidence tables, one can check that $\A'_1$, $\A_2$ and $\A'_2$ have the same intersection lattices as the one-parameter families $\mathcal{O}_{276}$, $\mathcal{O}_{155}$ and $\mathcal{O}_{154}$ from \cite{Me, CKC}, respectively. On the other hand, the intersection lattice of $\A_1$ has a two-dimensional realization space, so the double octics associated with $\A_1$ have $h^{1,2}\ge2$, and they do not appear in \cite{CKC}.

It is worth looking at the first pair more closely. Both arrangements have $8$ quadruple points, and in both cases every plane passes through exactly $4$ of them. The difference lies in the way these points are distributed. In $\A_1$ the quadruple points come in four pairs $\{P,P'\}$ such that the four planes through $P$ and the four planes through $P'$ give a partition of $\A_1$. For instance, the planes $x,y,z,x+y+z$ pass through $(0:0:0:1)$, and the remaining four planes pass through $(1:2:2:0)$. In $\A'_1$ there is no such pair of quadruple points.
\end{example}

\begin{table}[ht]
\centering
\small
\caption{Weak Ziegler pairs of eight planes. The last column gives the dimension of the realization space $\mathcal{R}$ of the intersection lattice, i.e., the number of moduli.}\label{tab:weak8}
\begin{tabular}{llllllc}
\toprule
 & $W$ & $\mathbf{d}$ & $\mathbf{c}$ & $\mathbf{b}$ & $t$ & $\dim\mathcal{R}$\\
\midrule
$\mathcal{O}_{241}$ & $(8;\,28;\,16,10)$ & $(3,4_6)$ & $(5_3,6_2)$ & $(7)$ & $4$ & $0$\\
$\mathcal{O}_{242}$ & $(8;\,28;\,16,10)$ & $(3,4_6)$ & $(5_3,6_2)$ & $(7)$ & $4$ & $1$\\
$\mathcal{O}_{239}$ & $(8;\,28;\,16,10)$ & $(4_{10})$ & $(5_9)$ & $(6_2)$ & $5$ & $0$\\
$\mathcal{O}_{240}$ & $(8;\,28;\,16,10)$ & $(4_{10})$ & $(5_9)$ & $(6_2)$ & $5$ & $0$\\
\midrule
$\A_1$ & $(8;\,28;\,24,8)$ & $(3,4_4,5_4)$ & $(5,6_8)$ & $(7_3)$ & $4$ & $2$\\
$\A'_1$ & $(8;\,28;\,24,8)$ & $(4_8)$ & $(5_3,6_4)$ & $(7_2)$ & $5$ & $1$\\
\midrule
$\A_2$ & $(8;\,25,1;\,16,11)$ & $(3,4_7,5)$ & $(5_6,6_2)$ & $(6,7)$ & $4$ & $1$\\
$\A'_2$ & $(8;\,25,1;\,16,11)$ & $(4_{11})$ & $(5_{11})$ & $(6_3)$ & $5$ & $1$\\
\bottomrule
\end{tabular}
\end{table}

\begin{remark}\label{rem:types}
Table~\ref{tab:weak8} shows that in each weak Ziegler pair of eight planes the arrangement with $\operatorname{mdr}=3$ has type $4$, and the arrangement with $\operatorname{mdr}=4$ has type $5$. This is in accordance with \cite[Theorem 5.1]{JP}, where it is shown that the types of rigid Cynk--Szemberg octic arrangements belong to $\{1,2,3,4,5\}$. Let us also recall that by Corollary~\ref{cor:octic} no Cynk--Szemberg octic arrangement has type $0$.
\end{remark}

\begin{example}\label{ex:weak9}
For $d=9$ let us present two weak Ziegler pairs. The first one consists of arrangements having only double lines:
\begin{align*}
\A_{3}:&\ f=xyz(x+y+z)w(4w-x-2y)(3x+2z+4w)(z-3y+4w)(z+w-4x-y),\\
\A'_{3}:&\ f=xyz(x+y+z)w(27x-y+w)(z+3w-8x)(z-62x-6y-2w)(3y+z+w-62x).
\end{align*}
Here $W(\A_3)=W(\A'_3)=(9;\,36;\,52,8)$, and
\[
\mathbf{d}(\A_3)=(4,5_4,6_{11}),\qquad \mathbf{d}(\A'_3)=(5_8,6_4),
\]
so $\operatorname{mdr}(\A_3)=4$ and $\operatorname{mdr}(\A'_3)=5$. The second one consists of arrangements having three triple lines:
\begin{align*}
\A_{4}:&\quad f=xy(x+y)z(3z-x)\,w\,(y+3z-3w)(y+3z+4w)(y+z+w-3x),\\
\A'_{4}:&\quad f=xy(x+y)z(z-3x)\,w\,(y-2z-4w)(y-2z+2w)(3x+2y+2z+2w).
\end{align*}
Here $W(\A_4)=W(\A'_4)=(9;\,27,3;\,25,10,3)$, and
\[
\mathbf{d}(\A_4)=(3,4_2,5_4,6),\qquad \mathbf{d}(\A'_4)=(4_5,5_4),
\]
so $\operatorname{mdr}(\A_4)=3$ and $\operatorname{mdr}(\A'_4)=4$. In Table~\ref{tab:weak9} we list the weak combinatorics of all $22$ weak Ziegler pairs of $9$ planes that we verified. In all these cases $t_4=t_5=0$ and $p_6=0$.
\end{example}

\begin{table}[ht]
\centering
\small
\caption{Weak combinatorics of verified weak Ziegler pairs of $9$ planes, together with the two values of $\operatorname{mdr}$.}\label{tab:weak9}
\begin{tabular}{cccccc}
\toprule
$t_2$ & $t_3$ & $p_3$ & $p_4$ & $p_5$ & $\operatorname{mdr}$\\
\midrule
$36$ & $0$ & $52$ & $8$ & $0$ & $4$ and $5$ \\
$36$ & $0$ & $48$ & $9$ & $0$ & $4$ and $5$ \\
$36$ & $0$ & $28$ & $14$ & $0$ & $4$ and $5$ \\
$33$ & $1$ & $41$ & $12$ & $0$ & $4$ and $5$ \\
$33$ & $1$ & $37$ & $13$ & $0$ & $4$ and $5$ \\
$33$ & $1$ & $29$ & $15$ & $0$ & $4$ and $5$ \\
$33$ & $1$ & $33$ & $14$ & $0$ & $4$ and $5$ \\
$27$ & $3$ & $25$ & $10$ & $3$ & $3$ and $4$ \\
$27$ & $3$ & $21$ & $11$ & $3$ & $3$ and $4$ \\
$27$ & $3$ & $17$ & $12$ & $3$ & $3$ and $4$ \\
$27$ & $3$ & $22$ & $8$ & $4$ & $3$ and $4$ \\
$27$ & $3$ & $18$ & $9$ & $4$ & $3$ and $4$ \\
$27$ & $3$ & $14$ & $10$ & $4$ & $3$ and $4$ \\
$27$ & $3$ & $19$ & $6$ & $5$ & $3$ and $4$ \\
$24$ & $4$ & $14$ & $8$ & $5$ & $3$ and $4$ \\
$24$ & $4$ & $20$ & $9$ & $4$ & $3$ and $4$ \\
$24$ & $4$ & $16$ & $10$ & $4$ & $3$ and $4$ \\
$24$ & $4$ & $12$ & $11$ & $4$ & $3$ and $4$ \\
$24$ & $4$ & $17$ & $7$ & $5$ & $3$ and $4$ \\
$24$ & $4$ & $13$ & $8$ & $5$ & $3$ and $4$ \\
$24$ & $4$ & $18$ & $4$ & $6$ & $3$ and $4$ \\
$24$ & $4$ & $16$ & $4$ & $6$ & $3$ and $4$ \\
\bottomrule
\end{tabular}
\end{table}

\begin{remark}\label{rem:noZiegler}
We would like to point out that we tried to find Ziegler pairs of essential and irreducible arrangements of at most nine planes which do not satisfy the condition (MDR), but we failed. For a given intersection lattice we computed $\operatorname{mdr}$ for many realizations of this lattice, chosen along random one-parameter families in the space of all realizations, and we checked whether $\operatorname{mdr}$ jumps somewhere. For $d\le8$ we did this for all intersection lattices, and for $d=9$ for about $90\%$ of them. The value of $\operatorname{mdr}$ never changed within a fixed intersection lattice. Of course, this is not a proof that such Ziegler pairs do not exist, but it strongly suggests that this is the case. Let us also stress that we compared only the values of $\operatorname{mdr}$, so our computations say nothing about Ziegler pairs with at most nine planes which satisfy the condition (MDR).

Let us give an easy explanation of this phenomenon. Let $P\in\PP^3_{\C}$ be a point of multiplicity $m_P$, i.e., exactly $m_P$ planes of $\A$ pass through $P$, and let $v\in\C^4$ be a vector representing $P$. Then the derivation
\[
\theta_P=\Big(\prod_{H\in\A,\ P\notin H}\ell_H\Big)\,\partial_v
\]
belongs to $D(\A)$ and it is not a multiple of $\theta_E$. Hence $\operatorname{mdr}(\A)\le d-m_P$. Combining this with Remark~\ref{rem:mdr12}, for every essential and irreducible arrangement of $d$ planes we obtain
\[
2\le\operatorname{mdr}(\A)\le d-m(\A),
\]
where $m(\A)$ denotes the maximal multiplicity of a point of $\A$. If $\A$ is not generic, that is, if some four planes of $\A$ pass through one point, then $m(\A)\ge4$, and for $d\le9$ we get $2\le\operatorname{mdr}(\A)\le5$. In other words, derivations which come directly from the combinatorics of $\A$ already live in very low degrees, and this leaves very little room for an additional derivation of even smaller degree, which would be caused by a special position of the planes rather than by their combinatorics. In Ziegler's example the situation is different: the arrangement has only triple points, so the above bound gives $\operatorname{mdr}\le9-3=6$, and the special position of the six vertices of the hexagon on a conic produces a derivation of degree $5$, just below this bound. On the other hand, for generic arrangements of planes the module $D(\A)$ depends only on $d$ \cite{RT, Y}, so generic arrangements cannot give Ziegler pairs.
\end{remark}

\begin{example}\label{ex:DP}
The mechanism described in Remark~\ref{rem:noZiegler} can be nicely seen in the Ziegler pair of Dimca and Pokora \cite[Example 4.1, Remark 4.5]{DP}. It consists of two arrangements of $11$ planes
\[
\A:f=w(w+L)g=0\qquad\text{and}\qquad\A':f'=w(w+L)g'=0,
\]
where $L=x-3y+5z$ and $g,g'\in\C[x,y,z]$ define two arrangements of $9$ lines with Ziegler's combinatorics, namely the six sides and the three main diagonals of a hexagon, such that the vertices of the hexagon lie on a conic for $g=0$, but not for $g'=0$. The arrangements $\A$ and $\A'$ are essential, irreducible, they have isomorphic intersection lattices, and
\[
\mathbf{d}(\A)=(2,6_2,7_7),\qquad \mathbf{d}(\A')=(2,7_{13}),
\]
see \cite[Example 4.1]{DP}. Hence $(\A,\A')$ is a Ziegler pair which satisfies the condition (MDR). The reason why $\operatorname{mdr}(\A)=\operatorname{mdr}(\A')=2$ is purely combinatorial. Indeed, the nine planes defined by the factors of $g$, resp. $g'$, pass through the point $O=(0:0:0:1)$, so $O$ is a $9$-fold point, and the bound from Remark~\ref{rem:noZiegler} gives $\operatorname{mdr}\le 11-9=2$. In fact, the derivation of degree $2$ found in \cite[Lemma 4.3]{DP} is
\[
\rho=(2w+L)E_0-w(10L+9w)\partial_w=(2w+L)\,\theta_E-11\,\theta_O,
\]
where $E_0=x\partial_x+y\partial_y+z\partial_z$ and $\theta_O=w(w+L)\partial_w$. On the other hand, the conic condition gives a Jacobian relation of degree $5$ of $g$, and by \cite[Lemma 4.3]{DP} this relation produces two Jacobian relations of degree $6$ of $f$, which are responsible for the entry $6_2$ in $\mathbf{d}(\A)$. So the special position of the hexagon is still visible in the module $D_0(\A)$, but only in degree $6$, which is far above the value of $\operatorname{mdr}$ forced by the $9$-fold point.
\end{example}

This leads us to the following natural problem.

\begin{problem}\leavevmode
\begin{enumerate}
\item Find a Ziegler pair of essential and irreducible plane arrangements in $\PP^3_{\C}$ which does not satisfy the condition (MDR). What is the minimal number of planes in such a pair?
\item Is there a Ziegler pair of essential and irreducible arrangements of at most $10$ planes in $\PP^3_{\C}$?
\end{enumerate}
\end{problem}

\section*{Acknowledgments}
We used the database of matroids from \cite{MMIB} and \textsc{Singular} \cite{Sing} to perform the computations described above. The \textsc{Singular} scripts verifying all computational claims of this note are available from the author upon request.


\begin{thebibliography}{99}

\bibitem{Abe}
T. Abe, Divisionally free arrangements of hyperplanes. \emph{Invent. Math.} \textbf{204}: 317--346 (2016).

\bibitem{AS}
C. Anghel, R. Stan, Terao's freeness conjecture for arrangements of 15 and 16 lines. \textbf{arXiv:2610.25868} (2026).

\bibitem{BK}
M. Barakat, L. K\"uhne, Computing the nonfree locus of the moduli space of arrangements and Terao’s freeness conjecture. \emph{Math. Comput.} \textbf{92(341)}: 1431 -- 1452 (2023).

\bibitem{ADP}
T. Abe, A. Dimca, P. Pokora, A new hierarchy for complex plane curves. \emph{Canad. Math. Bull.} \textbf{69(3)}: 824--848 (2026).

\bibitem{CKC}
S. Cynk, B. Kocel-Cynk, Classification of double octic Calabi--Yau threefolds with $h^{1,2}\le1$ defined by an arrangement of eight planes. \emph{Commun. Contemp. Math.} \textbf{22(1)}: 1850082, 38 pp. (2020).

\bibitem{CS}
S. Cynk, T. Szemberg, Double covers and Calabi--Yau varieties. In: \emph{Singularities Symposium -- {\L}ojasiewicz 70}, Banach Center Publ. \textbf{44}, Polish Acad. Sci. Inst. Math., Warsaw, 1998, 93--101.

\bibitem{DIM}
A. Dimca, D. Ibadula, A. M\u{a}cinic, Freeness for $13$ lines arrangements is combinatorial. \emph{Discrete Math.} \textbf{342} (2019).

\bibitem{DP}
A. Dimca, P. Pokora, On the Jacobian algebras of Ziegler pairs of plane arrangements. Preprint (2026).

\bibitem{DS}
A. Dimca, G. Sticlaru, From Pascal's Theorem to the geometry of Ziegler's line arrangements. \emph{J. Algebraic Combin.} \textbf{60}: 991--1009 (2024).

\bibitem{FV}
D. Faenzi, J. Vall\`es, Logarithmic bundles and line arrangements, an approach via the standard construction. \emph{J. London Math. Soc.} \textbf{90} (2014).

\bibitem{JP}
M. Janasz, P. Pokora, On homological properties of some Cynk--Szemberg octic hyperplane arrangements. \emph{J. Korean Math. Soc.} \textbf{63(4)}: 753--771 (2026).

\bibitem{Me}
C. Meyer, \emph{Modular Calabi--Yau threefolds}. Fields Institute Monographs \textbf{22}, Amer. Math. Soc., Providence, RI, 2005.

\bibitem{MMIB}
Y. Matsumoto, S. Moriyama, H. Imai, D. Bremner, Matroid enumeration for incidence geometry. \emph{Discrete Comput. Geom.} \textbf{47}: 17--43 (2012).

\bibitem{OT}
P. Orlik, H. Terao, \emph{Arrangements of hyperplanes}. Grundlehren der Mathematischen Wissenschaften \textbf{300}. Springer-Verlag, Berlin, 1992.

\bibitem{RT}
L. Rose, H. Terao, A free resolution of the module of logarithmic forms of a generic arrangement. \emph{J. Algebra} \textbf{136}: 376--400 (1991).

\bibitem{Sing}
W. Decker, G.-M. Greuel, G. Pfister, H. Sch\"onemann, \textsc{Singular} 4-3-2 --- A computer algebra system for polynomial computations. \url{https://www.singular.uni-kl.de} (2023).

\bibitem{Te}
H. Terao, Generalized exponents of a free arrangement of hyperplanes and Shepherd-Todd-Brieskorn formula. \emph{Invent. Math.} \textbf{63}: 159--179 (1981).

\bibitem{Y}
S. Yuzvinsky, A free resolution of the module of derivations for generic arrangements. \emph{J. Algebra} \textbf{136}: 432--438 (1991).

\bibitem{Z}
G. M. Ziegler, Multiarrangements of hyperplanes and their freeness. In: \emph{Singularities (Iowa City, IA, 1986)}, Contemp. Math. \textbf{90}, Amer. Math. Soc., Providence, RI, 1989, 345--359.

\bibitem{Z2}
G. M. Ziegler, Combinatorial construction of logarithmic differential forms. \emph{Adv. Math.} \textbf{76}: 116--154 (1989).

\end{thebibliography}
\end{document}